\documentclass{amsart}
\usepackage{amsxtra,amssymb,multicol,graphicx}
\begin{document}
\title{What are \dots\ multiple orthogonal polynomials?}
\author[A.~Mart\'\i nez-Finkelshtein and W.~Van Assche]{Andrei Mart\'\i nez-Finkelshtein and Walter Van Assche \\
Universidad de Almer\'\i a and KU Leuven}
\date{\today}


\rightline{\Huge \bf WHAT IS \dots\ A MULTIPLE}\bigskip
\rightline{\Huge \bf ORTHOGONAL POLYNOMIAL?}
\vskip1cm
\rightline{\it \huge Andrei Mart\'\i nez-Finkelshtein and Walter Van Assche}
\vskip1in

\begin{multicols}{2}

\textit{Multiple orthogonal polynomials} (MOPs) are polynomials of one variable which satisfy orthogonality conditions with respect to several measures. They should not be confused with multivariate of multivariable orthogonal polynomials, which are polynomials of several variables.  Other terminology is also used, e.g.,
Hermite-Pad\'e polynomials, polyorthogonal polynomials (Nikishin and Sorokin \cite{NikiSor}), $d$-orthogonal polynomials (the latter primarily by the French-Tunisian school of Maroni, Doauk, Ben Cheikh and collaborators).
They are a very useful extension of orthogonal polynomials, and recently received renewed interest because tools have become available to investigate  their asymptotic behavior  and they do appear in a number of
interesting applications. Various families of special MOPs have been found, extending the classical
orthogonal polynomials but also giving completely new special functions \cite{Aptekarev}, \cite[Ch.~23]{Ismail}. 

\begin{figure*}[t]\label{fig:zeros}
	\center \includegraphics[scale=1]{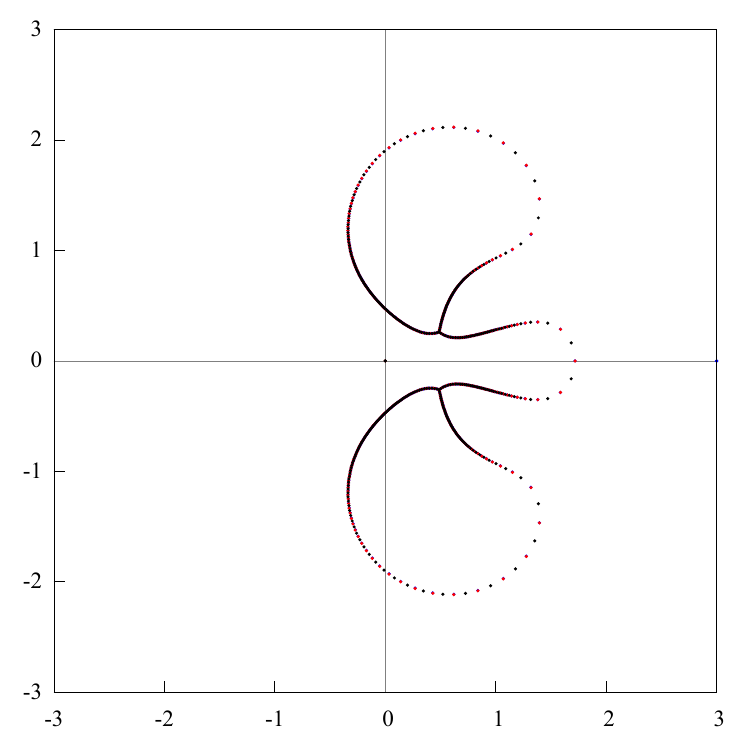}%
	\caption{Zeros of  $A_{\vec{n},1}$, $A_{\vec{n},2}$ and $B_{\vec{n} }$, $\vec{n}=(350, 350)$, corresponding to Type I Hermite-Pad\'e approximation to the functions $f_1=w$, $f_2=w^2$, where $w(z)$ is a solution of the algebraic equation $w^3+3(z-3)^2 w - 2i (3z-1)^3=0$ (picture courtesy of S.P.~Suetin).}
\end{figure*} 

\section{Definition}Let $r$ be a positive integer and $\mu_1,\mu_2,\ldots,\mu_r$ be positive measures on the real line for which all the moments
are finite. Let $\vec{n} = (n_1,n_2,\ldots,n_r) \in \mathbb{N}^r$ be a multi-index of size $|\vec{n}| = n_1+n_2+\cdots+n_r$.
There are two types of multiple orthogonal polynomials. \textit{Type I multiple orthogonal polynomials} are given as a vector
$(A_{\vec{n},1},A_{\vec{n},2},\ldots,A_{\vec{n},r})$ of $r$ polynomials, where $A_{\vec{n},j}$ has degree $\leq n_j-1$, for
which
\begin{equation}  \label{typeI}
   \sum_{j=1}^r \int x^k A_{\vec{n},j}(x)\, d\mu_j(x) = 0, \qquad 0 \leq k \leq |\vec{n}|-2, 
\end{equation}
and one usually adds the normalization
\begin{equation}  \label{typeInorm}
  \sum_{j=1}^r \int x^{|\vec{n}|-1}  A_{\vec{n},j}(x)\, d\mu_j(x) = 1 .   
\end{equation}
If all the measures are absolutely continuous with respect to one measure $\mu$, i.e.,  $d\mu_j(x) =w_j(x)\,d\mu(x)$, then
the function
\[   Q_{\vec{n}}(x) = \sum_{j=1}^r A_{\vec{n},j}(x)w_j(x)  \]
is orthogonal to all polynomials of degree $\leq |\vec{n}|-2$ with respect to the measure $\mu$.
The \textit{type II multiple orthogonal polynomial} $P_{\vec{n}}$ is the monic polynomial of degree $|\vec{n}|$ that satisfies
the orthogonality conditions
\begin{equation}  \label{typeII}    
\int P_{\vec{n}}(x) x^k \, d\mu_j(x) = 0, \qquad  0 \leq k \leq n_j-1, 
\end{equation}
for $1 \leq j \leq r$. 
Both conditions \eqref{typeI}--\eqref{typeInorm} and  conditions \eqref{typeII} yield a corresponding  linear system of $|\vec{n}|$ equations in the $|\vec{n}|$ unknowns, either
coefficients of the polynomials $(A_{\vec{n},1},\dots,A_{\vec{n},r})$ or coefficients of the monic polynomial $P_{\vec{n}}$.
The matrices of these linear systems are each other's transpose and contain moments of the $r$ measures $(\mu_1,\ldots,\mu_r)$.
A solution of these linear systems may not exist or may not be unique. One needs extra assumptions on the measures
$(\mu_1,\ldots,\mu_r)$ for a solution to exist and be unique. If a unique solution exists for a multi-index $\vec{n}$ then
the multi-index is said to be normal. If all multi-indices are normal, then the system $(\mu_1,\ldots,\mu_r)$ is said to be a perfect system.

One can also define MOPs of mixed type, combining both types of orthogonality conditions \eqref{typeI} and \eqref{typeII}.

\section{Hermite-Pad\'e approximation}
Multiple orthogonal polynomials originate from \textit{Hermite-Pad\'e approximation}, a 
simultaneous rational approximation to several functions.
Suppose $f_1,f_2,\ldots,f_r$ are analytic functions in a neighborhood of infinity, with  series expansions of the form
\[   f_j(z) = \sum_{k=0}^\infty \frac{c_{k,j}}{z^{k+1}}, \qquad 1 \leq j \leq r.  \]
Type I Hermite-Pad\'e approximation is to find polynomials $(A_{\vec{n},1},\ldots,A_{\vec{n},r})$ and $B_{\vec{n}}$,
with $\deg A_{\vec{n},j} \leq n_j-1$, such that
\[   \sum_{j=1}^r A_{\vec{n},j}(z)f_j(z) - B_{\vec{n}}(z) = \mathcal{O}\Bigl(\frac{1}{z^{|\vec{n}|}} \Bigr), \quad z \to \infty. \]
Type II Hermite-Pad\'e approximation is to find rational approximants $Q_{\vec{n},j}/P_{\vec{n}}$ with a common denominator 
for the functions $f_1,\ldots,f_r$ in the sense
that
\[    P_{\vec{n}}(z) f_j(z) - Q_{\vec{n},j}(z) = \mathcal{O}\Bigl( \frac{1}{z^{n_j+1}} \Bigr), \quad z\to \infty,  \]
for $1 \leq j \leq r$. If
\[   f_j(z) = \int \frac{d\mu_j(x)}{z-x}, \]
then $(A_{\vec{n},1},\ldots,A_{\vec{n},r})$ are the type I multiple orthogonal polynomials, and $P_{\vec{n}}$ is the type II
multiple orthogonal polynomial for the measures $(\mu_1,\ldots,\mu_r)$. 

In the case when $r=1$, both approximations coincide, and the rational function  $Q_{\vec{n},1}/P_{\vec{n}}$ is known simply as the \textit{Pad\'e approximant} to $f_1$ at $z=\infty$.

\section{Number theory}

The construction of the previous section goes back to the end of the 19th century, especially, to Hermite (and his student Pad\'e) and to Klein (and his scientific descendants, Lindemann and Perron). Hermite's proof that $e$ is a transcendental number is based on Hermite-Pad\'e approximation.
Using these ideas, Lindemann generalized this result proving that $e^{\alpha_1 z}, \dots ,e^{\alpha_r z}$ are algebraically independent over $\mathbb Q$  as long as the algebraic numbers $\alpha_1, \dots, \alpha_r$ are linearly independent over $\mathbb Q$ (which, in turn, yields  transcendence of $e$ and $\pi$).

Ap\'ery's more recent proof from 1979 that $\zeta(3)$ is irrational can be seen as a problem of Hermite-Pad\'e approximation of the functions
\[   f_1(z) = \int_0^1 \frac{1}{z-x}\, dx, \quad f_2(z) = - \int_0^1 \frac{\log x}{z-x}\, dx, \]
\[  f_3(z) = \frac12 \int_0^1 \frac{\log^2 x}{z-x}\, dx, \]
in the sense that one wants to find polynomials $A_n$, $B_n$, $C_n$, $D_n$ of degree $\leq n$ for which (as $z \to \infty$)
\begin{align*}
    A_n(z) f_1(z) - B_n(z) f_2(z) - C_n(z) &= \mathcal{O}\left( \frac{1}{z^{n+1}} \right), \\
    A_n(z) f_2(z) - 2B_n(z) f_3(z) - D_n(z) &= \mathcal{O}\left( \frac{1}{z^{n+1}} \right),
\end{align*}
with the extra condition $A_n(1)=0$.
The polynomials $(A_n,B_n)$ then are a vector of mixed type I--type II multiple orthogonal polynomials, which can be obtained explicitly
in terms of Legendre polynomials on $[0,1]$. The relevant quantity is $f_3(1)=\zeta(3)$, and the construction gives rational approximants
to $\zeta(3)$ which are better than possible if $\zeta(3)$ would be rational, implying that this number is irrational.

Similar Hermite-Pad\'e constructions were used by Ball and Rivoal in 2001 to show that infinitely many $\zeta(2n+1)$ are irrational,
and somewhat later Zudilin was able to prove that at least one of the numbers $\zeta(5), \zeta(7), \zeta(9), \zeta(11)$ is irrational \cite{Fischler}.
Krattenthaler, Rivoal and Zudilin also obtained irrationality results for a $q$-extension of the zeta function in 2006, and the rational approximation
now involves multiple little $q$-Jacobi polynomials.  

\section{Random matrices}
A research field where multiple orthogonal polynomials have appeared more recently and turned out to be very useful, is random matrix theory.
It was well known that orthogonal polynomials play an important role in the so-called Gaussian Unitary Ensemble (GUE) of random matrices.
These are random $n\times n$ hermitian matrices $M$ for which the probability distribution of the matrix entries is given by a density
\[       \frac{1}{Z_n} e^{-n \textup{Tr} V(M)} \, dM, \]
where $V$ is of polynomial growth at infinity and $Z_n$ is a normalizing constant that makes this a probability density. The average characteristic polynomial 
\[    P_n(z) = \mathbb{E}\det (zI_n -M)  \]
is a monic polynomial of degree $n$ that satisfies
\[   \int P_n(x) x^k e^{-V(x)}\, dx = 0, \qquad k=0,1,\ldots,n-1, \]
and hence it is the $n$th degree orthogonal polynomial for the weight $e^{-V(x)}$.
Another random matrix model is one with an external source \cite{Kuijlaars}, i.e., a fixed (non-random) hermitian matrix $A$, and the
probability density now is
\[   \frac{1}{Z_n} e^{-n \textup{Tr}(M^2-AM)}\, dM. \]
If $A$ has $n_j$ eigenvalues $a_j$ $(1 \leq j \leq r)$ and $n=|\vec{n}|$, then the average characteristic polynomial
\[    P_{\vec{n}}(z) = \mathbb{E}\det (zI_n - M)  \]
is a type II multiple orthogonal polynomial with orthogonality conditions
\[   \int_{-\infty}^{\infty} P_{\vec{n}}(x) x^k e^{-n(x^2-a_jx)}\, dx = 0, \quad 0 \leq k \leq n_j-1, \]
for $1 \leq j \leq r$, and this is a so-called \textit{multiple Hermite polynomial} \cite{Delvaux}. 
Instead of Hermitian matrices, one can also use positive definite matrices from the Wishart ensemble and find multiple Laguerre polynomials
\cite{Bleher}.
More recently, products of Ginibre random matrices were investigated by Akemann, Ipsen and Kieburg in 2013, and the singular values of such matrices
are described in terms of multiple orthogonal polynomials for which the weight functions are Meijer G-functions \cite{Kuijlaars-Zhang}.

\section{Non-intersecting random paths}

\begin{figure*}[t]
	\center \includegraphics[scale=0.45]{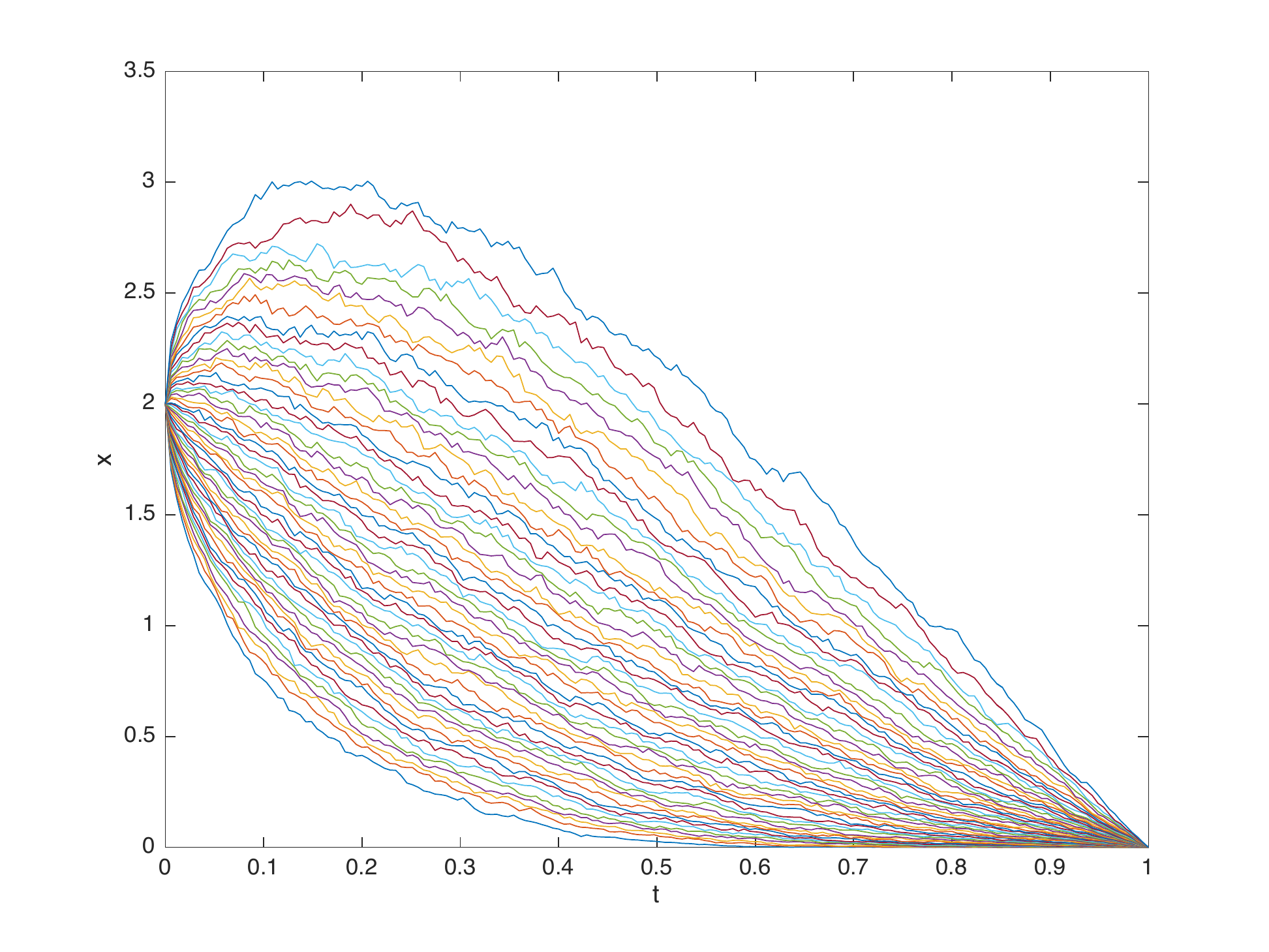}%
	\caption{Numerical simulation of 50 rescaled non-intersecting BESQ, starting at $a = 2$ and ending at the origin. The statistics of the position of the paths at each time $t$ is described in terms of MOPs with respect to weights related to modified Bessel functions $I_\nu$ and $I_{\nu-1}$.}
	\label{fig:BESQ}
\end{figure*} 

Eigenvalues and singular values of random matrices are special cases of determinantal point processes. These are
point processes for which the correlation function can be written as a determinant of a kernel function $K$:
\[    \rho_n(x_1,\ldots,x_n) = \det \bigl( K(x_i,x_j) \bigr)_{1 \leq i,j \leq n}  \]
for every $n\geq 1$. For the eigenvalues of a random matrix from  GUE the kernel function is
\[   K_n(x,y) = e^{-n(V(x)+V(y))/2} \sum_{k=0}^{n-1} p_k(x)p_k(y), \]
where the $(p_k)_k$ are the orthonormal polynomials for the weight $e^{-V(x)}$ on $\mathbb{R}$. This kernel is the well-known Christoffel-Darboux kernel for orthogonal polynomials.  For random matrices with an external source the kernel is in terms of type I and type II multiple orthogonal polynomials
\[    K_n(x,y) = \sum_{k=0}^{n-1}  P_{\vec{n}_k}(x)Q_{\vec{n}_{k+1}}(y), \]
where $(\vec{n}_k)_{0\leq k \leq n}$ is a path in $\mathbb{N}^r$ from $\vec{0}$ to $\vec{n}$, and the sum turns out to be independent of this path.

Another determinantal point process is related to non-intersecting random paths, in particular non-intersecting Brownian motions \cite{Daems} 
and squared Bessel paths (BESQ) \cite{Kuijlaars-MF-Wielonsky}. As Kuijlaars and collaborators have shown, non-intersecting Brownian motions leaving at $r$ points and arriving
at 1 point can be described in terms of multiple Hermite polynomials and the analysis uses the same tools as random matrices with an external source.
If the Brownian motions leave from $r$ points and arrive at $r$ points, then one needs a mixture of type I and type II multiple orthogonal polynomials.
BESQ are Brownian motions in $d$-dimensional space, in which one tracks the distance of the particle from the origin. This distance squared is distributed as a non-central $\chi^2$ random variable, and the dependence structure (repulsion of the particles, since the paths may not intersect) gives a determinantal structure with multiple orthogonal polynomials related to modified Bessel functions $I_\nu$ and $I_{\nu-1}$, where $\nu=d/2$, see Figure~\ref{fig:BESQ}.

\section{Integrable systems}
Multiple orthogonal polynomials satisfy a system of \textit{nearest neighbor recurrence relations}, extending the well-known three term recurrence relation for orthogonal polynomials. For this reason, MOPs are related to the spectral theory of non-symmetric operators (see the work of Aptekarev, Kalyagin and collaborators). For the same reason, MOPs are also appearing in the theory of integrable systems \cite{Aptekarev-D-VA}. To mention a few examples, multiple Jacobi polynomials
were used by Mukhin and Varchenko \cite{Mukhin} to give a counterexample to the Bethe Ansatz conjecture for the Gaudin model, 
some non-Hermitian oscillator Hamiltonians were constructed using multiple Charlier \cite{Miki2} and multiple Meixner polynomials \cite{Miki1}, 
multiple Meixner-Pollaczek polynomials appear in the six-vertex model \cite{Bender},
and discrete integrable systems can be generated by multiple orthogonal polynomials.

\section{Analytic theory of MOPs}

Since the classical families of orthogonal polynomials (Jacobi, Laguerre, Hermite) satisfy so many properties, there are several available tools to study their analytic behavior (zeros, asymptotics, etc.), which are absent when considering other families. The first ``universal'' tool in this sense was \textit{logarithmic potential theory}: it is known (this knowledge can be traced all the way back to Gauss, but was fully developed at the end of 20th century in the works of Saff, Mhaskar, Gonchar, Rakhmanov, Stahl and others) that the global zero (or weak) asymptotics of general families of orthogonal polynomials can be described in terms of solutions of extremal problems for the logarithmic energy of positive measures (with multiple variations, such as the presence of an external field, upper constraints, symmetry conditions, etc.). In many cases, these solutions are related to harmonic functions on compact \textit{Riemann surfaces}. For a more detailed description of the asymptotics, especially at some particular points on the plane, we have to recur to special functions (both classical and modern, such as Painlev\'e transcendents). 

Recently a new method, developed by Deift and Zhou, based on the \textit{matrix Riemann-Hilbert characterization} of orthogonal polynomials (due to Fokas, Its, and Kitaev) has emerged, which encompasses these previous tools and, when successful, gives detailed information about the asymptotic behavior, uniformly in the whole complex plane.

As the pioneering work of Nikishin, followed by Aptekarev, Gonchar, Kalyagin, Kuijlaars, Rakhmanov, Suetin, Van Assche and others showed, there are direct extensions of all these techniques for the study of multiple orthogonal polynomials: one now needs to consider extremal problems for the logarithmic energy for a vector of measures \cite[\S 5.4]{NikiSor}, the related Riemann surfaces now typically have a larger number of sheets, and the matrix Riemann-Hilbert characterization involves matrices of higher order \cite{VA-Ger-Kuijl}. These are not merely technical complications: the asymptotic behavior of multiple orthogonal polynomials is so rich that it is difficult to envision any general theory in this respect, at least in the near future.

\end{multicols}

\noindent Andrei Mart\'\i nez-Finkelshtein \\
Universidad de Almer\'\i a, Spain \\
andrei@ual.es
\medskip

\noindent Walter Van Assche \\
KU Leuven, Belgium \\
walter@wis.kuleuven.be 

\end{document}